\documentclass[10pt]{amsart}
\usepackage{amsfonts}
\usepackage{amsmath}

\newtheorem{theorem}{Theorem}
\newtheorem{definition}{Definition}
\newtheorem{proposition}{Proposition}

\begin{document}
\title{Symmetries and conservation laws of Hamiltonian systems}
\author{Liviu Popescu}
\maketitle

\begin{abstract}
In this paper we study the infinitesimal symmetries, Newtonoid vector
fields, infinitesimal Noether symmetries and conservation laws of
Hamiltonian systems. Using the dynamical covariant derivative and Jacobi
endomorphism on the cotangent bundle we find the invariant equations of
infinitesimal symmetries and Newtonoid vector fields and prove that the
canonical nonlinear connection induced by a regular Hamiltonian can be
determined by these symmetries. Finally, an example from optimal control
theory is given.
\end{abstract}



\section{\textbf{Introduction}}

The notion of symmetry plays a very important role in all field theories
being related with conservation laws by Noether type theorems. The use of
the symmetries of a system in the description of its dynamical evolution has
a long history and goes back to the classical mechanics (see for instance
\cite{Ab},\cite{Ar},\cite{Mrs}). Also, the Lagrangian and Hamiltonian
formalisms are fundamental concepts in physics, differential equations or
optimal control and in most of cases the study starts with a variational
problem formulated for a regular Lagrangian on the tangent bundle $TM$ over
the manifold $M$ and very often the whole set of problems is transferred on
the dual space $T^{*}M$,$\,$ endowed with a Hamiltonian function, via
Legendre transformation.

The present paper contains some contributions to the study of
symmetries of Hamiltonian systems and shows how the well-known
local symmetries of Lagrangian systems emerge in Hamiltonian
formulation. The tangent bundle has a canonical tangent structure
$J$ and together with a semispray $S$ (system of second order
differential equation-SODE) induce a nonlinear connection that
describes the geometry of the system \cite{Cr1},\cite{Gr1}. These
structures lead to the notions of Jacobi endomorphism and
dynamical covariant derivative introduced by J. Carine\~na and E.
Mart\'\i nez (see \cite{Ca},\cite{Ma}) which have been used in the
study of symmetries for SODE in \cite{Bu3}. But, the existence of
a symplectic structure on the tangent bundle depends on a
Lagrangian function on $TM$. We have to remark that some type of
symmetries on the tangent bundle as dynamical symmetries, Lie
symmetries, Cartan (Noether) symmetries (see e.g.
\cite{Bu1},\cite{Cr3},\cite{Le},\cite{Le1},\cite{Pr1},\cite{Pr2},\cite{Sa1})
and Newtonoid vector field (\cite{Bu3},\cite{Mar}) have been
studied in a lot of papers, where the main geometric structures
are the semispray, the symplectic structure $\omega _{L}$induced
by a regular Lagrangian $L$ and the energy $E_L$. Contrary, the
cotangent bundle is endowed with a canonical symplectic structure
and does not have a canonical tangent structure or something
similar with a semispray. However, the existence of a
pseudo-metric structure or a regular Hamiltonian on $T^{*}M$
permit us to define an adapted tangent structure $\mathcal{J}$ and
a regular vector field $\rho $ which induce a nonlinear connection
\cite{Op}. These geometrical structures permit us to introduce the
Jacobi endomorphism and dynamical covariant derivative on $T^{*}M$
(see \cite{Po2},\cite{Po3}) which will be used in this paper in
order to find the invariant equations of the infinitesimal
symmetries of Hamiltonian systems. In fact, this work containts
the ideas proposed by the author in \cite{Po3}. Different types of
symmetries and conservation laws for Hamiltonian systems can be
found, for example, in \cite
{Ba},\cite{Cas},\cite{Ja},\cite{Mrl},\cite{Mu},\cite{Pu}.

The paper is organized as follows. In the second section the
preliminary geometric structures on the cotangent bundle are
recalled (see for instance \cite{Mi2},\cite{Op},\cite{Po1},
\cite{Po2},\cite{Po3},\cite{Ya} and references therein). We
introduce the Berwald linear connection $\mathcal{D}$ on $T^{*}M$
induced by a nonlinear connection $\mathcal{N}$ and study its
properties in subsection 2.1. We show that this linear connection
is compatible with the horizontal and vertical projectors, adapted
tangent structure and complex structure. Also, we find its action
on the local Berwald basis. Moreover, we prove that in the case of
the horizontal $ \mathcal{J}$-regular vector field $\rho $, the
Berwald linear connection coincides with the dynamical covariant
derivative, that is $\mathcal{D}_\rho =\nabla $. Consequently, in
this case, the integral curves of a $\mathcal{J}$-regular vector
field are the geodesics of the Berwald linear connection.

In the third section we investigate the symmetries of Hamiltonian systems on
the cotangent bundle. First, for a regular Hamiltonian on $T^{*}M,$ we
introduce an integrable adapted tangent structure $\mathcal{J}_H$, a regular
vector field $\rho _H$, which is the Hamiltonian vector field, and find the
coefficients of the canonical nonlinear connection. Moreover, we give the
expression of the Jacobi endomorphism, which depends only on the regular
Hamiltonian and find the action of the dynamical covariant derivative on the
local Berwald basis. Next, using the Hamiltonian vector field, canonical
symplectic structure and adapted tangent structure, we study the
infinitesimal symmetries, natural infinitesimal symmetries, Newtonoid vector
field, infinitesimal Noether symmetries and conservation laws of Hamiltonian
systems. Also, using the dynamical covariant derivative and Jacobi
endomorphism on $T^{*}M$, we find the invariant equations of the
infinitesimal symmetries and Newtonoid vector field and prove that these
symmetries determine the canonical nonlinear connection. Moreover, we show
when one of these symmetries will imply the others and that there is a one
to one correspondence between the exact infinitesimal Noether symmetry and
conservation laws. Finally, an example from optimal control theory is given.

\section{\textbf{Geometrical structures on the cotangent bundle}}

Let $M$ be a differentiable, $n$-dimensional manifold and $(T^{*}M,\tau ,M)$
its cotangent bundle. If the local coordinates on $\tau ^{-1}(U)$ are
denoted $(x^i,p_i),$ $(i=\overline{1,n})$ then the natural basis on $T^{*}M$
is $\left( \frac \partial {\partial x^i},\frac \partial {\partial
p_i}\right) $ and $(dx^i,dp_i)$ is the dual natural basis. The following
geometric objects \cite{Ya}
\begin{equation*}
C^{*}=p_i\frac \partial {\partial p_i},\quad \theta =p_idx^i,\quad \omega
=d\theta =dp_i\wedge dx^i,
\end{equation*}
have the following properties:

1$^{\circ }$ $C^{*}$ is a vertical vector field, globally defined
on $T^{*}M$ , which is called the Liouville-Hamilton vector field.

2$^{\circ }$ The 1-form $\theta $ is globally defined on $T^{*}M$ and is
called the Liouville 1-form.

3$^{\circ }$ The 2-form $\omega $ is the canonical symplectic structure.\\If
$L$ and $K$ are $(1,1)$-type tensor field then the Fr\"olicher-Nijenhuis
bracket $[L,K]$ is the vector valued 2-form \cite{Fr}
\begin{eqnarray*}
\lbrack L,K](X,Y) &=&[LX,KY]+[KX,LY]+(LK+KL)[X,Y] \\
&&\ \ \ \ -L[X,KY]-K[X,LY]-L[KX,Y]-K[LX,Y].
\end{eqnarray*}
and the Nijenhuis tensor of $L$ is given by
\begin{equation*}
N_L(X,Y)=\frac 12[L,L]=[LX,LY]+L^2[X,Y]-L[X,LY]-L[LX,Y].
\end{equation*}
For a vector field $X$ in $\mathcal{X}(M)$ the
Fr\"olicher-Nijenhuis bracket $[X,L]=\mathcal{L}_XL$ is the
$(1,1)$-type tensor field on $M$ given by $
\mathcal{L}_XL=\mathcal{L}_X\circ L-L\circ \mathcal{L}_X$, where
$\mathcal{L} _X$ is the usual Lie derivative. On the cotangent
bundle $T^{*}M$ there exists the integrable vertical distribution
$V_uT^{*}M$, $u\in T^{*}M$ generated locally by the basis $\left\{
\frac \partial {\partial p_i}\right\} _{i=\overline{1,n}}$. A
nonlinear connection on $T^{*}M$ is defined by an almost product
structure $\mathcal{N}$ (i.e. a morphism $
\mathcal{N}:\mathcal{X}(T^{*}M)\rightarrow \mathcal{X}(T^{*}M)$
with $ \mathcal{N}^2=Id$) such that $VT^{*}M=Ker\left(
Id+\mathcal{N}\right) $. If $ \mathcal{N}$ is a nonlinear
connection then $HT^{*}M=Ker(Id-\mathcal{N})$ is the horizontal
distribution associated to $\mathcal{N}$, which is supplementary
to the vertical distribution, that is $TT^{*}M=VT^{*}M\oplus
HT^{*}M.$ If $\mathcal{N}$ is a nonlinear connection then on the
every domain of the local chart $\tau ^{-1}(U)$, the adapted basis
of the horizontal distribution $HT^{*}M$ is \cite{Mi2}
\begin{equation*}
\frac \delta {\delta x^i}=\frac \partial {\partial
x^i}+\mathcal{N} _{ij}\frac \partial {\partial p_j},
\end{equation*}
where $\mathcal{N}_{ij}(x,p)$ are the coefficients of the nonlinear
connection. The dual adapted basis is $\delta p_i=dp_i-\mathcal{N}_{ij}dx^j.$
The system of vector fields $\left( \frac \delta {\delta x^i},\frac \partial
{\partial p_i}\right) $ defines the local Berwald basis on $T^{*}M$. A
nonlinear connection induces the horizontal and vertical projectors given by
\begin{equation*}
h=\frac 12(Id+\mathcal{N}),\quad v=\frac 12(Id-\mathcal{N}),\quad h=\frac
\delta {\delta x^i}\otimes dx^i,\quad v=\frac \partial {\partial p_i}\otimes
\delta p_i,
\end{equation*}
which satisfy the properties $h^2=h$, $v^2=v$, $hv=vh=0$,
$h+v=Id$, $h-v= \mathcal{N}$. The nonlinear connection
$\mathcal{N}$ on $T^{*}M$ is called symmetric if $\omega
(hX,hY)=0,$ for $X,Y\in \mathcal{X}(T^{*}M)$, that is $
\mathcal{N}_{ij}=\mathcal{N}_{ji}.$ The following equations hold
\cite{Mi2}
\begin{equation}
\left[ \frac \delta {\delta x^i},\frac \delta {\delta x^j}\right]
=R_{ijk}\frac \partial {\partial p_k},\quad \ \left[ \frac \delta
{\delta x^i},\frac \partial {\partial p_j}\right] =-\frac{\partial
\mathcal{N}_{ir}}{
\partial p_j}\frac \partial {\partial p_r},\quad \left[ \frac \partial
{\partial p_i},\frac \partial {\partial p_j}\right] =0,
\end{equation}
\begin{equation}
R_{ijk}=\frac{\delta \mathcal{N}_{jk}}{\delta x^i}-\frac{\delta
\mathcal{N} _{ik}}{\delta x^j}.
\end{equation}
The curvature of the nonlinear connection $\mathcal{N}$ on
$T^{*}M$ is given by $\Omega =-\frac 12[h,h],$ where $h$ is the
horizontal projector and $ \frac 12[h,h]$ is the Nijenhuis tensor
of $h$. In local coordinates $\Omega =-\frac 12R_{ijk}dx^i\wedge
dx^j\otimes \frac \partial {\partial p_k},$ where $R_{ijk}$ is
given by (2). An almost tangent structure on $T^{*}M$ is a
morphism $\mathcal{J}:\mathcal{X}(T^{*}M)$ $\rightarrow
\mathcal{X} (T^{*}M) $ of rank $n$ such that $\mathcal{J}^2=0$.
The almost tangent structure is called adapted if
$Im\mathcal{J}=Ker\mathcal{J}=VT^{*}M$ (see \cite{Op}). The
following properties hold
\begin{equation}
\mathcal{J}h=\mathcal{J},\quad h\mathcal{J}=0,\quad
\mathcal{J}v=0,\quad v \mathcal{J}=\mathcal{J}.
\end{equation}
Locally, an adapted almost tangent structure has the form
\begin{equation}
\mathcal{J}=t_{ij}dx^i\otimes \frac \partial {\partial p_j},
\end{equation}
where $t_{ij}(x,p)$ is a tensor field of rank $n.$ The existence
of a nonlinear connection on $T^{*}M$ is equivalent with the
conditions $\mathcal{ NJ}=-\mathcal{J}$,
$\mathcal{JN}=\mathcal{J}$. The adapted almost tangent structure
$\mathcal{J}$ is integrable if and only if $\frac{\partial
t^{ij}}{
\partial p_k}=\frac{\partial t^{kj}}{\partial p_i}$, where $
t_{ij}t^{jk}=\delta _i^k.$ Also, $\mathcal{J}$ is called symmetric
if $ \omega (\mathcal{J}X,Y)=\omega (\mathcal{J}Y,X)$, which
locally is equivalent with the symmetry of the tensor
$t_{ij}(x,p)$. From \cite{Op}, \cite{Po3} we have that a vector
field $\rho $ $\in \mathcal{X}(T^{*}M)$ is called
$\mathcal{J}$-regular if it satisfies the equation
\begin{equation}
\mathcal{J}[\rho ,\mathcal{J}X]=-\mathcal{J}X,\quad \forall X\in
\mathcal{X} (T^{*}M).
\end{equation}
Locally, a vector field on $T^{*}M$ given in local coordinates by
\begin{equation*}
\rho =\xi ^i(x,p)\frac \partial {\partial x^i}+\chi _i(x,p)\frac \partial
{\partial p_i},
\end{equation*}
is $\mathcal{J}$-regular if and only if $t^{ij}=\frac{\partial \xi
^j}{
\partial p_i}$, where $t_{ij}t^{jk}=\delta _i^k.$ For a $\mathcal{J}$
-regular vector field $\rho $ and an arbitrary nonlinear
connection $ \mathcal{N}$ with induced $(h,v)\,$ projectors, we
consider the vertically valued $(1,1)$-type tensor field on
$T^{*}M\backslash \{0\}$ given by \cite {Po3}
\begin{equation*}
\Phi =v\circ \mathcal{L}_\rho h,
\end{equation*}
which is called the Jacobi endomorphism. In local coordinates we obtain
\begin{equation*}
\mathcal{L}_\rho \frac \delta {\delta x^j}=-\frac{\delta \xi
^i}{\delta x^j} \frac \delta {\delta x^i}+R_{jk}\frac \partial
{\partial p_k},\quad \mathcal{ L}_\rho \frac \partial {\partial
p_j}=-t^{ji}\frac \delta {\delta x^i}+\left(
t^{ji}\mathcal{N}_{ik}-\frac{\partial \chi _k}{\partial p_j}
\right) \frac \partial {\partial p_k}.
\end{equation*}
Locally, the Jacobi endomorphism has the form
\begin{equation}
\Phi =\mathcal{R}_{ij}dx^i\otimes \frac \partial {\partial
p_j},\quad \mathcal{R}_{jk}=\frac{\delta \xi ^i}{\delta
x^j}\mathcal{N}_{ik}-\frac{ \delta \chi _k}{\delta x^j}+\rho
(\mathcal{N}_{jk}).
\end{equation}
We can also recover the Jacobi endomorphism $\Phi $ from the
curvature tensor $\Omega $ through the formula $\Phi =i_\rho
\Omega +v\circ \mathcal{L} _{v\rho }h.$ Moreover, if $\rho $ is a
horizontal $\mathcal{J}$-regular vector field then $\rho =h\rho $,
$v\rho =0$ and $\Phi =i_\rho \Omega $. Locally, it results
\cite{Po3}
\begin{equation*}
\rho =\xi ^i\frac \delta {\delta x^i},\quad \chi _i=\xi ^k\mathcal{N}%
_{ki},\quad \mathcal{R}_{ij}=R_{kij}\xi ^k,
\end{equation*}
which show us the relation between the Jacobi endomorphism given
by (6) and curvature tensor from (2). For a given
$\mathcal{J}$-regular vector field $ \rho $ on $T^{*}M$ the Lie
derivative $\mathcal{L}_\rho $ defines a tensor derivation on
$T^{*}M$, but does not preserve some of the geometric structure as
adapted tangent structure or nonlinear connection. Next, using a
nonlinear connection, we introduce a tensor derivation on
$T^{*}M$, called the dynamical covariant derivative, that
preserves some geometric structures (see e.g.
\cite{Bu3},\cite{Ma},\cite{Sa2} for the tangent bundle case).
Using \cite{Po3} we set:

\begin{definition}
A map $\nabla :\mathcal{T}(T^{*}M\backslash \{0\})\rightarrow
\mathcal{T} (T^{*}M\backslash \{0\})$ is said to be a tensor
derivation on $ T^{*}M\backslash \{0\}$ if the following
conditions are satisfied:

i) $\nabla $ is $\Bbb{R}$-linear,

ii) $\nabla $ is type preserving, i.e. $\nabla (\mathcal{T}%
_s^r(T^{*}M\backslash \{0\}))\subset
\mathcal{T}_s^r(T^{*}M\backslash \{0\})$ , for each $(r,s)\in
\Bbb{N}\times \Bbb{N.},$

iii) $\nabla $ obeys the Leibnitz rule $\nabla (T\otimes S)=\nabla T\otimes
S+T\otimes \nabla S,$

iv) $\nabla \,$commutes with any contractions.
\end{definition}

For a $\mathcal{J}$-regular vector field $\rho $ and an arbitrary
nonlinear connection $\mathcal{N}$ with induced $(h,v)\,$
projectors, we consider the map $\nabla
:\mathcal{X}(T^{*}M\backslash \{0\})\rightarrow \mathcal{X}
(T^{*}M\backslash \{0\})$ given by
\begin{equation}
\nabla =h\circ \mathcal{L}_\rho \circ h+v\circ \mathcal{L}_\rho \circ v,
\end{equation}
which is called the dynamical covariant derivative with respect to $\rho $
and the nonlinear connection $\mathcal{N}$. By setting $\nabla f=\rho (f),$
for $f\in C^\infty (T^{*}M\backslash \{0\})$ using the Leibnitz rule and the
requirement that $\nabla $ commutes with any contraction, we can extend the
action of $\nabla $ to arbitrary tensor fields on $T^{*}M\backslash \{0\}$
(see \cite{Po3}). By direct computation we obtain $\nabla h=\nabla v=0$ and
the action of $\nabla $ on the Berwald basis:
\begin{equation}
\nabla \frac \delta {\delta x^j}=-\frac{\delta \xi ^i}{\delta x^j}\frac
\delta {\delta x^i},\quad \nabla dx^j=\frac{\delta \xi ^j}{\delta x^i}dx^i.
\end{equation}
\begin{equation}
\nabla \frac \partial {\partial p_j}=\left(
t^{ji}\mathcal{N}_{ik}-\frac{
\partial \chi _k}{\partial p_j}\right) \frac \partial {\partial p_k},\quad
\nabla \delta p_j=-\left( t^{ki}\mathcal{N}_{ij}-\frac{\partial
\chi _j}{
\partial p_k}\right) \delta p_k.
\end{equation}
The following results hold \cite{Po3}
\begin{equation}
h\circ \mathcal{L}_\rho \circ \mathcal{J}=-h,\quad
\mathcal{J}\circ \mathcal{ L}_\rho \circ v=-v,
\end{equation}
\begin{equation}
\nabla \mathcal{J}=\mathcal{L}_\rho \mathcal{J}+h-v,\quad \nabla
\mathcal{J} =\left( \rho (t_{ij})+t_{kj}\frac{\partial \xi
^k}{\partial x^i}-t_{ik}\frac{
\partial \chi _j}{\partial p_k}+2\mathcal{N}_{ij}\right) dx^i\otimes \frac
\partial {\partial p_j}.
\end{equation}
Given an adapted tangent structure $\mathcal{J}$ and a $\mathcal{J}$-regular
vector field $\rho $, then the compatibility condition $\nabla \mathcal{J}=0$
fix the canonical nonlinear connection with $h$, $v$ projectors
\begin{equation*}
h=\frac 12\left( Id-\mathcal{L}_\rho \mathcal{J}\right) ,\quad v=\frac
12\left( Id+\mathcal{L}_\rho \mathcal{J}\right) ,
\end{equation*}
The (1,1)-type tensor field
\begin{equation}
\mathcal{N}=-\mathcal{L}_\rho \mathcal{J},
\end{equation}
is the almost product structure which will be used in the following. The
local coefficients are given by \cite{Op}
\begin{equation}
\mathcal{N}_{ij}=\frac 12\left( t_{ik}\frac{\partial \chi
_j}{\partial p_k} -t_{kj}\frac{\partial \xi ^k}{\partial x^i}-\rho
(t_{ij})\right) .
\end{equation}
The almost complex structure has the form $\Bbb{F}=h\circ
\mathcal{L}_\rho h- \mathcal{J}$ and in local coordinates we have
$\Bbb{F}=t^{ij}\frac \delta {\delta x^i}\otimes \delta
p_j-t_{ij}\frac \partial {\partial p_i}\otimes dx^j$. The
dynamical covariant derivative has in this case the properties
\cite{Po3}
\begin{equation*}
\nabla \mathcal{J}=0,\ \nabla \Bbb{F}=0.
\end{equation*}
Moreover, if $\rho $ is a horizontal $\mathcal{J}$-regular vector field then
$\nabla \rho =0$.

\subsection{\textbf{Berwald linear connection on the cotangent bundle}}

Next, we introduce the Berwald linear connection induced by a nonlinear
connection and prove that in the case of horizontal $\mathcal{J}$-regular
vector field it coincides with the dynamical covariant derivative. This
connection was introduced on the tangent bundle by L. Berwald in \cite{Be}
and studied later in \cite{Cr2}, \cite{Mi2} and \cite{Bu2}.

\begin{definition}
The Berwald linear connection on the cotangent bundle is given by
\begin{equation*}
\mathcal{D}:\mathcal{X}(T^{*}M\backslash \{0\})\times \mathcal{X}
(T^{*}M\backslash \{0\})\rightarrow \mathcal{X}(T^{*}M\backslash
\{0\}),
\end{equation*}
\begin{equation}
\mathcal{D}_XY=v[hX,vY]+h[vX,hY]+\mathcal{J}[vX,(\Bbb{F}+\mathcal{J})Y]+(
\Bbb{F}+\mathcal{J})[hX,\mathcal{J}Y].
\end{equation}
\end{definition}

Because all the structures from the right hand side of (14) are
additive, it results that $\mathcal{D}$ is also additive, with
respect to both arguments. Next, we prove that $D_{fX}Y=fD_XY$,
$\forall f\in C^\infty (T^{*}M)$ by straightforward computation,
using the relations $vh=hv=\mathcal{J}v=0$ and (
$\Bbb{F}+\mathcal{J})h=0$. Indeed, for the first term from
$D_{fX}Y$ we have $v[fhX,vY]=v\left( f[hX,vY]-(vY)(f)hX\right) $
$=fv[hX,vY]$. In order to
prove the relation $D_XfY=X(f)Y+fD_XY$ we remark that $%
D_XfY=fD_XY+(hX)(f)v^2Y+(vX)(f)h^2Y+(vX)(f)\mathcal{J}(\Bbb{F}+\mathcal{J}
)(Y)+(hX)(f)(\Bbb{F}+\mathcal{J})\mathcal{J}(Y).$ But $v^2=v$,
$h^2=h$, $ \mathcal{J}(\Bbb{F}+\mathcal{J})=v$,
($\Bbb{F}+\mathcal{J})\mathcal{J}=h$
(see \cite{Po3}) and it results $%
D_XfY=fD_XY+(hX)(f)vY+(vX)(f)hY+(vX)(f)vY+(hX)(f)hY=fD_XY+(hX)(f)Y+(vX)(f)Y=fD_XY+X(f)Y
$ which prove that $\mathcal{D}$ is a linear connection.

\begin{proposition}
The Berwald linear connection has the following properties
\begin{equation*}
\mathcal{D}h=0,\quad \mathcal{D}v=0,\quad \mathcal{DJ}=0,\quad
\mathcal{D} \Bbb{F}=0.
\end{equation*}
\end{proposition}

\textbf{Proof}. Using the properties of the vertical and
horizontal projectors we obtain
$\mathcal{D}_XvY=v[hX,vY]+\mathcal{J}[vX,(\Bbb{F}+ \mathcal{J})Y]$
and $v(\mathcal{D}_XY)=v[hX,vY]+\mathcal{J}[vX,(\Bbb{F}+
\mathcal{J})Y]$ which yields $\mathcal{D}v=0$. Also, $\mathcal{D}
_XhY=h[vX,hY]+(\Bbb{F}+\mathcal{J})[hX,JY]=h(\mathcal{D}_XY)$ and
it results $\mathcal{D}h=0$. \\Moreover,
$\mathcal{D}_X\mathcal{J}Y=v[hX,\mathcal{J}Y]+ \mathcal{J}[vX,hY]$
and $\mathcal{J}(\mathcal{D}_XY)=\mathcal{J}[vX,hY]+v[hX,
\mathcal{J}Y]$ and we obtain $\mathcal{D}\mathcal{J}=0.$
$\mathcal{D}_X\Bbb{F
}Y=v[hX,-\mathcal{J}Y]+h[vX,(\Bbb{F}+\mathcal{J})Y]+\mathcal{J}[vX,-hY]+(
\Bbb{F}+\mathcal{J})[hX,vY]$ and
$\Bbb{F}(\mathcal{D}_XY)=(\Bbb{F}+\mathcal{J
})[hX,vY]-\mathcal{J}[vX,hY]+h[vX,(\Bbb{F}+\mathcal{J})Y]-v[hX,\mathcal{J}Y]=
\mathcal{D}_X\Bbb{F}Y$ which yields $\mathcal{D}\Bbb{F}=0.$ \hfill
\hbox{\rlap{$\sqcap$}$\sqcup$}\\It results that the Berwald
connection preserves both horizontal and vertical vector fields.
Locally, we have the following formulas
\begin{eqnarray*}
\mathcal{D}_{\frac \delta {\delta x^i}}\frac \delta {\delta x^j}
&=&t^{ks}\left( \frac{\delta t_{jk}}{\delta
x^i}-t_{jr}\frac{\partial \mathcal{N}_{ik}}{\partial p_r}\right)
\frac \delta {\delta x^s},\quad \mathcal{D}_{\frac \delta {\delta
x^i}}\frac \partial {\partial p_j}=-\frac{
\partial \mathcal{N}_{ir}}{\partial p_j}\frac \partial {\partial p_r},\quad
\\
\mathcal{D}_{\frac \partial {\partial p_i}}\frac \delta {\delta x^j}
&=&0,\quad \mathcal{D}_{\frac \partial {\partial p_i}}\frac \partial
{\partial p_j}=t_{ks}\frac{\partial t^{jk}}{\partial p_i}\frac \partial
{\partial p_s}.
\end{eqnarray*}
We can see that the dynamical covariant derivative has the same properties
and this leads to the next result.

\begin{theorem}
If $\rho $ is a horizontal $\mathcal{J}$-regular vector field then the
following equality holds
\begin{equation*}
\nabla =\mathcal{D}_\rho .
\end{equation*}
\end{theorem}

\textbf{Proof}. If $\rho $ is a horizontal $\mathcal{J}$-regular vector
field then $\rho =h\rho $ and $v\rho =0$ which implies
\begin{equation*}
\mathcal{D}_\rho Y=v[\rho ,vY]+(\Bbb{F}+\mathcal{J})[\rho ,\mathcal{J}Y].
\end{equation*}
But $\nabla Y=h[\rho ,hY]+v[\rho ,vY]$ and we will prove that
$h[\rho ,hY]=( \Bbb{F}+\mathcal{J})[\rho ,\mathcal{J}Y]$ using the
computation in local coordinates. Let us consider $Y=X^i(x,p)\frac
\partial {\partial x^i}+Y_j(x,p)\frac \partial {\partial p_j}$ and
using (1) we get
\begin{equation*}
\lbrack \rho ,hY]=\left[ \xi ^i\frac \delta {\delta x^i},X^j\frac
\delta {\delta x^j}\right] =\xi ^iX^jR_{ijk}\frac \partial
{\partial p_k}+\xi ^i \frac{\delta X^j}{\delta x^i}\frac \delta
{\delta x^j}-X^j\frac{\delta \xi ^i }{\delta x^j}\frac \delta
{\delta x^i},
\end{equation*}

\begin{equation*}
h[\rho ,hY]=\left( \xi ^i\frac{\delta X^j}{\delta
x^i}-X^i\frac{\delta \xi ^j }{\delta x^i}\right) \frac \delta
{\delta x^j}.
\end{equation*}
Introducing the expression of the canonical nonlinear connection in the case
of horizontal $\mathcal{J}$-regular vector field given by
\begin{equation*}
\mathcal{N}_{ij}=-t_{kj}\frac{\partial \xi ^k}{\partial x^i}-\xi
^k\frac{ \delta t_{ij}}{\delta x^k}+\xi ^lt_{ik}\frac{\partial
\mathcal{N}_{lj}}{
\partial p_k},
\end{equation*}
we obtain
\begin{eqnarray*}
h[\rho ,hY] &=&\left( \xi ^i\frac{\delta X^j}{\delta
x^i}-X^i\left( \frac{
\partial \xi ^j}{\partial x^i}+\left( -t_{sk}\frac{\partial \xi ^s}{\partial
x^i}-\xi ^s\frac{\delta t_{ik}}{\delta x^s}+\xi ^lt_{is}\frac{\partial
\mathcal{N}_{lk}}{\partial p_s}\right) t^{kj}\right) \right) \frac \delta
{\delta x^j} \\
\ &=&\left( \xi ^i\frac{\delta X^j}{\delta x^i}-X^i\left(
\frac{\partial \xi ^j}{\partial x^i}-\frac{\partial \xi
^j}{\partial x^i}-\xi ^st^{kj}\frac{ \delta t_{ik}}{\delta
x^s}+\xi ^lt_{is}\frac{\partial \mathcal{N}_{lk}}{
\partial p_s}t^{kj}\right) \right) \frac \delta {\delta x^j} \\
\ &=&\left( \xi ^i\frac{\delta X^j}{\delta x^i}+X^i\xi
^s\frac{\delta t_{ik} }{\delta x^s}t^{kj}-\xi
^lX^it_{is}\frac{\partial \mathcal{N}_{lk}}{\partial
p_s}t^{kj}\right) \frac \delta {\delta x^j}.
\end{eqnarray*}
Next
\begin{equation*}
\lbrack \rho ,\mathcal{J}Y]=\left[ \xi ^i\frac \delta {\delta
x^i},X^kt_{kl}\frac \partial {\partial p_l}\right] =-\xi
^iX^kt_{kl}\frac{
\partial \mathcal{N}_{ij}}{\partial p_l}\frac \partial {\partial p_j}+\xi
^i\frac \delta {\delta x^i}(X^kt_{kl})\frac \partial {\partial
p_l}-X^kt_{kl} \frac{\partial \xi ^i}{\partial p_l}\frac \delta
{\delta x^i},
\end{equation*}
and using that $(\Bbb{F}+\mathcal{J})\left( \frac \delta {\delta x^i}\right)
=0,$ $(\Bbb{F}+\mathcal{J})\left( \frac \partial {\partial p_j}\right)
=t^{jk}\frac \delta {\delta x^i}$ we obtain
\begin{equation*}
(\Bbb{F}+\mathcal{J})[\rho ,\mathcal{J}Y]=\left( \xi
^i\frac{\delta X^j}{ \delta x^i}+X^k\xi ^i\frac{\delta
t_{kl}}{\delta x^i}t^{lj}-\xi ^iX^kt_{kl} \frac{\partial
\mathcal{N}_{is}}{\partial p_l}t^{sj}\right) \frac \delta {\delta
x^j},
\end{equation*}
which ends the proof. \hfill
\hbox{\rlap{$\sqcap$}$\sqcup$}\\Moreover, $ \nabla \rho
=\mathcal{D}_\rho \rho =0$ and it results that the integral curves
of $\rho $ are geodesics of the Berwald linear connection.

\section{\textbf{Symmetries of Hamiltonian systems}}

A Hamilton space \cite{Mi2} is a pair $(M,H)$ where $M$ is a differentiable,
$n-$dimensional manifolds and $H$ is a function on $T^{*}M$ with the
properties:

1$^{\circ }$ $H:(x,p)\in T^{*}M\rightarrow H(x,p)\in \Bbb{R}$ is
differentiable on $T^{*}M$ and continue on the null section of the
projection $\tau :T^{*}M\rightarrow M$.

2$^{\circ }$ The Hessian of $H$ with respect to $p_i$ is nondegenerate
\begin{equation}
g^{ij}=\frac{\partial ^2H}{\partial p_i\partial p_j},\ \ rank\left\|
g^{ij}(x,p)\right\| =n\text{ on }T^{*}M\backslash \{0\}.
\end{equation}

3$^{\circ }$ The tensor field $g^{ij}(x,p)$ has constant signature
on $ T^{*}M\backslash \{0\}.$ \\The triple $(T^{*}M,\omega ,H)$ is
called a Hamiltonian system.

The Hamiltonian $H$ on $T^{*}M$ induces a pseudo-Riemannian metric $g_{ij}$
with $g_{ij}g^{jk}=\delta _i^k$ and $g^{jk}$ given by (15) on $VT^{*}M$. It
induces a unique adapted tangent structure denoted
\begin{equation*}
\mathcal{J}_H=g_{ij}dx^i\otimes \frac \partial {\partial p_j}.
\end{equation*}
A $\mathcal{J}$-regular vector field induced by the regular Hamiltonian $H$
has the form
\begin{equation*}
\rho _H=\frac{\partial H}{\partial p_i}\frac \partial {\partial x^i}+\chi
_i\frac \partial {\partial p_i}.
\end{equation*}
There exists a unique Hamiltonian vector field $\rho _H\in
\mathcal{X} (T^{*}M)$ which is a $\mathcal{J}$-regular vector
field such that $i_{\rho _H}\omega =-dH$, given by
\begin{equation}
\rho _H=\frac{\partial H}{\partial p_i}\frac \partial {\partial
x^i}-\frac{
\partial H}{\partial x^i}\frac \partial {\partial p_i}.
\end{equation}
The symmetric nonlinear connection
\begin{equation*}
\mathcal{N}=-\mathcal{L}_{\rho _H}\mathcal{J}_H,
\end{equation*}
has the coefficients (see \cite{Mi2}, \cite{Po1})
\begin{equation}
\mathcal{N}_{ij}=\frac 12\left( \{g_{ij},H\}-\left(
g_{ik}\frac{\partial ^2H }{\partial p_k\partial
x^j}+g_{jk}\frac{\partial ^2H}{\partial p_k\partial x^i}\right)
\right) ,
\end{equation}
where the Poisson bracket is
\begin{equation*}
\{g_{ij},H\}=\frac{\partial g_{ij}}{\partial p_k}\frac{\partial H}{\partial
x^k}-\frac{\partial H}{\partial p_k}\frac{\partial g_{ij}}{\partial x^k},
\end{equation*}
and is called the \textit{canonical nonlinear connection} of the Hamilton
space $(M,H),$ which is a metric nonlinear connection, that is $\nabla g=0$
(see \cite{Po2}). In this case, the coefficients of the Jacobi endomorphism
have the form
\begin{equation}
\mathcal{R}_{jk}=\frac{\partial ^2H}{\partial p_i\partial
x^j}\mathcal{N} _{ik}+\frac{\partial ^2H}{\partial p_i\partial
x^k}\mathcal{N}_{ji}+\mathcal{
N}_{jl}\mathcal{N}_{ik}g^{li}+\frac{\partial ^2H}{\partial
x^j\partial x^k} +\rho _H(\mathcal{N}_{jk}),
\end{equation}
and the action of the dynamical covariant derivative on the Berwald basis is
given by
\begin{equation}
\nabla \frac \delta {\delta x^j}=h\left[ \rho _H,\frac \delta
{\delta x^j}\right] =-\left( \frac{\partial ^2H}{\partial
p_i\partial x^j}+\mathcal{N }_{jk}g^{ki}\right) \frac \delta
{\delta x^i},
\end{equation}
\begin{equation}
\nabla \frac \partial {\partial p_j}=v\left[ \rho _H,\frac
\partial {\partial p_j}\right] =\left( \frac{\partial
^2H}{\partial p_j\partial x^i}+ \mathcal{N}_{ik}g^{kj}\right)
\frac \partial {\partial p_i}.
\end{equation}
In the following, we study the symmetries of Hamiltonian systems (see also
\cite{Pu}) on the cotangent bundle using the Hamiltonian vector field and
the adapted tangent structure.

\begin{definition}
A vector field $X\in \mathcal{X}(T^{*}M)$ is an infinitesimal symmetry of
Hamiltonian vector field if $[\rho _H,X]=0.$
\end{definition}

If we consider $X=X^i(x,p)\frac \partial {\partial x^i}+Y_i(x,p)\frac
\partial {\partial p_i}$ then an infinitesimal symmetry is given by the
equations
\begin{equation*}
X\left( \frac{\partial H}{\partial p_i}\right) =\rho _H\left( X^i\right)
,\quad \ X\left( \frac{\partial H}{\partial x^i}\right) +\rho _H\left(
Y^i\right) =0,
\end{equation*}
and the first relation leads to
\begin{equation*}
Y_k=g_{ki}\left( \rho _H\left( X^i\right) -X^j\frac{\partial ^2H}{\partial
p_i\partial x^j}\right) .
\end{equation*}

\begin{definition}
A vector field $\widetilde{Z}\in \mathcal{X}(M)$ is said to be a natural
infinitesimal symmetry if its complete lift to $T^{*}M$ is an infinitesimal
symmetry, that is $[\rho _H,\widetilde{Z}^{C*}]=0.$
\end{definition}

We know that, for $\widetilde{Z}=\widetilde{Z}^i(x)\frac \partial {\partial
x^i\text{ }}$ the complete lift on $T^{*}M$ is given by \cite{Ya}
\begin{equation*}
\widetilde{Z}^{C*}=\widetilde{Z}^i(x)\frac \partial {\partial
x^i}+p_j\frac{
\partial \widetilde{Z}^j}{\partial x^i}\frac \partial {\partial p_i},
\end{equation*}
and a natural infinitesimal symmetry is characterized by the equations
\begin{equation*}
\widetilde{Z}^{C*}\left( \frac{\partial H}{\partial p_k}\right)
=\frac{
\partial H}{\partial p_i}\frac{\partial \widetilde{Z}^k}{\partial x^i},
\end{equation*}
\begin{equation*}
\widetilde{Z}^{C*}\left( \frac{\partial H}{\partial x^k}\right)
=\frac{
\partial H}{\partial x^i}\frac{\partial \widetilde{Z}^i}{\partial x^k}-p_j
\frac{\partial H}{\partial p_i}\frac{\partial ^2\widetilde{Z}^j}{\partial
x^k\partial x^i}.
\end{equation*}
Next, we introduce the Newtonoid vector field on $T^{*}M$ (see
\cite{Mu}, \cite{Bu3} for tangent bundle case) which help us to
find the canonical nonlinear connection induced by a regular
Hamiltonian.

\begin{definition}
A vector field $X\in \mathcal{X}(T^{*}M)$ is called Newtonoid vector field
if $\mathcal{J}_H[\rho _H,X]=0.$
\end{definition}

In local coordinates we obtain
\begin{equation*}
g_{ij}\left( X\left( \frac{\partial H}{\partial p_i}\right) -\rho _H\left(
X^i\right) \right) \frac \partial {\partial p_j}=0,
\end{equation*}
and using that $\ rank\left\| g^{ij}(x,p)\right\| =n$ it result the equation
\begin{equation*}
X\left( \frac{\partial H}{\partial p_i}\right) =\rho _H\left( X^i\right) ,
\end{equation*}
which leads to the expression of a Newtonoid vector field
\begin{equation*}
X=X^i\frac \partial {\partial x^i}+g_{ki}\left( \rho _H\left(
X^i\right) -X^j \frac{\partial ^2H}{\partial p_i\partial
x^j}\right) \frac \partial {\partial p_k}.
\end{equation*}
We remark that $X$ is an infinitesimal symmetry if and only if it is
Newtonoid vector field and satisfies the equation
\begin{equation*}
X\left( \frac{\partial H}{\partial x^i}\right) +\rho _H\left( Y^i\right) =0.
\end{equation*}
The set of Newtonoid vector fields is given by
\begin{equation*}
\frak{X}_{\rho _H}=Ker\left( \mathcal{J}_H\circ \mathcal{L}_{\rho _H}\right)
=Im\left( Id+\mathcal{J}_H\circ \mathcal{L}_{\rho _H}\right) .
\end{equation*}
In the following, we will use the dynamical covariant derivative and Jacobi
endomorphism in order to find the invariant equations of Newtonoid vector
field and infinitesimal symmetries. Let $\rho _H$ be the Hamiltonian vector
field, $\mathcal{N}$ an arbitrary nonlinear connection with induced $v,h$
projectors and $\nabla $ the induced dynamical covariant derivative. We set:

\begin{proposition}
A vector field $X\in \mathcal{X}(T^{*}M)$ is a Newtonoid vector field if and
only if
\begin{equation}
v(X)=\mathcal{J}_H(\nabla X).
\end{equation}
\textbf{Proof}. We have the relation (10) $\mathcal{J}_H\circ
\nabla = \mathcal{J}_H\circ \mathcal{L}_{\rho _H}+v$ and it
results $\mathcal{J} _H[\rho _H,X]=0$ if and only if
$v(X)=\mathcal{J}_H(\nabla X).$ \hfill
\hbox{\rlap{$\sqcap$}$\sqcup$}
\end{proposition}

\begin{proposition}
A vector field $X\in \mathcal{X}(T^{*}M)$ is a infinitesimal symmetry if and
only if $X$ is a Newtonoid vector field and satisfies the equation
\begin{equation}
\nabla (\mathcal{J}_H\nabla X)+\Phi (X)=0.
\end{equation}
\end{proposition}

\textbf{Proof}. A vector field $X\in \mathcal{X}(T^{*}M)$ is a
infinitesimal symmetry if and only if $h[\rho _H,X]=0$ and $v[\rho
_H,X]=0$. Composing by $ \mathcal{J}_H$ we obtain
$\mathcal{J}_Hh[\rho _H,X]=\mathcal{J}_H[\rho _H,X]=0$ which means
that $X$ is a Newtonoid vector field. Also,
\begin{equation*}
v[\rho _H,X]=v[\rho _H,vX]+v[\rho _H,hX]=\nabla (vX)+\Phi
(X)=\nabla ( \mathcal{J}_H(\nabla X))+\Phi (X),
\end{equation*}
which ends the proof.\hfill
\hbox{\rlap{$\sqcap$}$\sqcup$}

For $f\in C^\infty (T^{*}M)$ and $X\in \mathcal{X}(T^{*}M)$ we define the
product
\begin{equation*}
f*X=\left( Id+\mathcal{J}_H\circ \mathcal{L}_{\rho _H}\right)
(fX)=fX+f \mathcal{J}_H[\rho _H,X]+\rho _H(f)\mathcal{J}_HX,
\end{equation*}
and it result that a vector field $X$ is a Newtonoid if and only if
\begin{equation*}
X=X^i(x,p)*\frac \partial {\partial x^i}.
\end{equation*}
Also, if $X\in \frak{X}_{\rho _H}$ then $f*X=fX+\rho _H(f)\mathcal{J}_HX$
(see \cite{Sa1},\cite{Bu3} for the case of tangent bundle). Next theorem
proves that the canonical nonlinear connection induced by a regular
Hamiltonian can be determined by symmetries.

\begin{theorem}
Let us consider the Hamiltonian vector field $\rho _H$, an arbitrary
nonlinear connection $\mathcal{N}$ and $\nabla $ the dynamical covariant
derivative. The following conditions are equivalent:

$i)$ $\nabla $ restricts to $\nabla :\frak{X}_{\rho _H}\rightarrow
\frak{X} _{\rho _H}$ satisfies the Leibnitz rule with respect to
the $*$ product.

$ii)$ $\nabla \mathcal{J}_H=0$

$iii)$ $\mathcal{L}_{\rho _H}\mathcal{J}_H+\mathcal{N}=0,$

$iv)$ $\mathcal{N}_{ij}=\frac 12\left( \{g_{ij},H\}-\left(
g_{ik}\frac{
\partial ^2H}{\partial p_k\partial x^j}+g_{jk}\frac{\partial ^2H}{\partial
p_k\partial x^i}\right) \right) .$
\end{theorem}

\textbf{Proof}. For $ii)\Rightarrow i)$ let us consider $X\in
\frak{X}_{\rho _H}$ and using (21) we get $vX=\mathcal{J}_H(\nabla
X)$. Applying $\nabla $ to both sides, we obtain $\nabla
(vX)=\nabla (\mathcal{J}_H\nabla X)$ which yields $(\nabla
v)X+v(\nabla X)=(\nabla \mathcal{J}_H)(\nabla X)+\mathcal{J}
_H\nabla (\nabla X)$. Using the relations $\nabla v=0$, $\nabla
\mathcal{J} _H=0$ it results $v(\nabla X)=\mathcal{J}_H\nabla
(\nabla X)$ which implies $ \nabla X\in \frak{X}_{\rho _H}$. For
$X\in \frak{X}_{\rho _H}$ we obtain
\begin{equation*}
\nabla \left( f*X\right) =\nabla (fX+\rho
_H(f)\mathcal{J}_HX)=\rho _H(f)X+f\nabla X+\rho
_H^2(f)\mathcal{J}_HX+\rho _H(f)\nabla (\mathcal{J} _HX),
\end{equation*}
\begin{equation*}
\nabla f*X+f*\nabla X=\rho _H(f)X+f\nabla X+\rho _H^2(f)\mathcal{J}_HX+\rho
_H(f)\mathcal{J}_H(\nabla X).
\end{equation*}
But $\nabla (\mathcal{J}_HX)=(\nabla
\mathcal{J}_H)X+\mathcal{J}_H(\nabla X)$ and from $\nabla
\mathcal{J}_H=0$ it results $\nabla (\mathcal{J}_HX)=
\mathcal{J}_H(\nabla X)$ which leads to $\nabla \left( f*X\right)
=\nabla f*X+f*\nabla X.$

For $i)\Rightarrow ii)$ we prove that $\nabla \mathcal{J}_H$
vanishes on the set $\frak{X}_{\rho _H}\cup
\mathcal{X}^v(T^{*}M\backslash \{0\})$ which is a set of
generators for $\mathcal{X}(T^{*}M\backslash \{0\})$. For $X\in
\mathcal{X}^v(T^{*}M\backslash \{0\})$ we have $\mathcal{J}_HX=0$
and $ \mathcal{J}_H(\nabla X)=0$ which lead to $\nabla
\mathcal{J}_H(X)=\nabla ( \mathcal{J}_HX)-\mathcal{J}_H(\nabla
X)=0$. Next, if $X\in \frak{X}_{\rho _H} $ then from $\nabla
\left( f*X\right) =\nabla f*X+f*\nabla X$ it results $\rho
_H(f)\nabla (\mathcal{J}_HX)=\rho _H(f)\mathcal{J}_H(\nabla X)$,
which implies $\rho _H(f)(\nabla \mathcal{J}_H)X=0$ for an
arbitrary function $ f\in C^\infty (T^{*}M\backslash \{0\})$ and
arbitrary vector field $X\in \frak{X}_{\rho _H}$. Therefore
$\nabla \mathcal{J}_H=0$ which ends the proof. The equivalence of
$ii)$, $iii)$, $iv)$ results from (11).\hfill
\hbox{\rlap{$\sqcap$}$\sqcup$}

Considering the canonical nonlinear connection
$\mathcal{N}=-\mathcal{L} _{\rho _H}\mathcal{J}_H$, we get the
following results.

\begin{proposition}
A vector field $X\in \mathcal{X}(T^{*}M)$ is a infinitesimal symmetry if and
only if $X$ is a Newtonoid vector field and satisfies the equation
\begin{equation}
\nabla ^2\mathcal{J}_HX+\Phi (X)=0,
\end{equation}
which locally yields
\begin{equation}
\nabla ^2g_{ij}X^i+\mathcal{R}_{ij}X^i=0.
\end{equation}
\end{proposition}

\textbf{Proof}. If $\mathcal{N}$ is the canonical nonlinear connection, then
$\nabla \mathcal{J}_H=0$ and using $\nabla \mathcal{J}_H=\nabla \circ
\mathcal{J}_H-\mathcal{J}_H\circ \nabla $ from (22) it results (23). Also,
we obtain that the local components of the vertical vector field (23) are
(24).\hfill
\hbox{\rlap{$\sqcap$}$\sqcup$}

\begin{definition}
a) An infinitesimal Noether symmetry of the Hamiltonian $H$ is a vector
field $X\in \mathcal{X}(T^{*}M)$ such that
\begin{equation*}
\mathcal{L}_X\omega =0,\quad \mathcal{L}_XH=0.
\end{equation*}

b) A vector field $\widetilde{X}\in \mathcal{X}(M)$ is said to be an
invariant vector field for the Hamiltonian $H$ if $\widetilde{X}^{C*}(H)=0.$

c) A function $f\in C^\infty (M)$ is a constant of motion (or a conservation
law) for the Hamiltonian $H$ if $\mathcal{L}_{\rho _H}f=0.$
\end{definition}

\begin{proposition}
Every infinitesimal Noether symmetry is an infinitesimal symmetry.
\end{proposition}

\textbf{Proof}. From the symplectic equation $i_{\rho _H}\omega =-dH,$
applying the Lie derivative in both sides, it results
\begin{equation*}
\mathcal{L}_X(i_{\rho _H}\omega )=-\mathcal{L}_XdH=-d\mathcal{L}_XH=0.
\end{equation*}
Also, from the formula $i_{[X,\rho _H]}=\mathcal{L}_X\circ i_{\rho
_H}-i_{\rho _H}\circ \mathcal{L}_X$ we obtain
\begin{equation*}
\mathcal{L}_X(i_{\rho _H}\omega )=i_{[X,\rho _H]}\omega +i_{\rho
_H}\mathcal{ L}_X\omega =i_{[X,\rho _H]}\omega ,
\end{equation*}
which leads to $i_{[X,\rho _H]}\omega =0$ and we get $[X,\rho _H]=0.$\hfill
\hbox{\rlap{$\sqcap$}$\sqcup$}

\begin{proposition}
If $\widetilde{X}$ is a vector field on $M$ such that
$\mathcal{L}_{ \widetilde{X}^{c*}}\theta _{\text{ }}$is closed and
$d(\widetilde{X} ^{c*}H)=0 $, then $\widetilde{X}$ is a natural
infinitesimal symmetry.
\end{proposition}

\textbf{Proof}. We have
\begin{eqnarray*}
i_{[\widetilde{X}^{c*},\rho _H]}\omega
&=&\mathcal{L}_{\widetilde{X} ^{c*}}(i_{\rho _H}\omega )-i_{\rho
_H}(\mathcal{L}_{\widetilde{X} ^{c*}}\omega
)=-\mathcal{L}_{\widetilde{X}^{c*}}dH-i_{\rho _H}(\mathcal{L}_{
\widetilde{X}^{c*}}d\theta )= \\
&=&-d\mathcal{L}_{\widetilde{X}^{c*}}H-i_{\rho
_H}d(\mathcal{L}_{\widetilde{X }^{c*}}\theta
)=-d(\widetilde{X}^{c*}H)=0,
\end{eqnarray*}
because $d(\mathcal{L}_{\widetilde{X}^{c*}}\theta )=0.$\hfill
\hbox{\rlap{$\sqcap$}$\sqcup$}

\begin{proposition}
The Hamiltonian vector field $\rho _H$ is an infinitesimal Noether symmetry.
\end{proposition}

\textbf{Proof}. Using the skew symmetry of the symplectic 2-form $\omega $,
it results
\begin{equation*}
0=i_{\rho _H}\omega (\rho _H)=-dH(\rho _H)=\rho _H(H)=\mathcal{L}_{\rho _H}H.
\end{equation*}
Also, from $d\omega =0$ we get
\begin{equation*}
\mathcal{L}_{\rho _H}\omega =di_{\rho _H}\omega +i_{\rho _H}d\omega
=-d(dH)=0.
\end{equation*}
\hfill
\hbox{\rlap{$\sqcap$}$\sqcup$}

Since Lie and exterior derivatives commute, we obtain for an infinitesimal
Noether symmetry
\begin{equation*}
d\mathcal{L}_X\theta =\mathcal{L}_Xd\theta =\mathcal{L}_X\omega =0.
\end{equation*}
It results that the 1-form $\mathcal{L}_X\theta $ is a closed 1-form and
consequently $\mathcal{L}_{\rho _H}\theta $ is closed.

\begin{definition}
An infinitesimal Noether symmetry $X$ $\in \mathcal{X}(T^{*}M)$ is said to
be an exact infinitesimal Noether symmetry if the 1-form $\mathcal{L}%
_X\theta $ is exact.
\end{definition}

The next result proves that there is a one to one correspondence between the
exact infinitesimal Noether symmetry and conservation laws. Also, if $X$ is
an exact infinitesimal Noether symmetry, then there is a function $f\in
C^\infty (M)$ such that $\mathcal{L}_X\theta =df$.

\begin{theorem}
If $X$ is an exact infinitesimal Noether symmetry, then $f-\theta (X)$ is a
conservation law for the Hamiltonian $H$. Conversely, if $f\in C^\infty (M)$
is a conservation law for $H$, then $X\in \mathcal{X}(T^{*}M\backslash
\{0\}) $ the unique solution of the equation $i_X\omega =-df$ is an exact
infinitesimal Noether symmetry.
\end{theorem}

\textbf{Proof}. We have $\rho _H(f-\theta (X))=d(f-\theta
(X))(\rho _H)=\left( \mathcal{L}_X\theta -di_X(\theta )\right)
(\rho _H)=i_Xd\theta (\rho _H)=i_X\omega (\rho _H)=-i_{\rho
_H}\omega (X)=dH(X)=0,$ and it results that $f-\theta (X)$ is a
conservation law for the dynamics associated to the regular
Hamiltonian $H$. Conversely, if $X$ is the solution of the
equation $i_X\omega =-df$ then $\mathcal{L}_X{\theta }=i_Xd{
\theta +}di_X\theta =-df+di_X\theta $ is an exact 1-form.
Consequently, $0=d \mathcal{L}_X\theta =\mathcal{L}_Xd\theta
=\mathcal{L}_X\omega .$ Also, $f$ is a conservation law, and we
have $0=\rho _H(f)=df(\rho _H)=-i_X\omega (\rho _H)=i_{\rho
_H}\omega (X)=-dH(X)=-X(H).$ Therefore, we obtain $
\mathcal{L}_XH=0$ and $X$ is an exact infinitesimal Noether
symmetry.\hfill \hbox{\rlap{$\sqcap$}$\sqcup$}

\begin{theorem}
If $\widetilde{X}\in \mathcal{X}(M)$ is an invariant vector field for the
Hamiltonian $H$ then its complete lift $\widetilde{X}^{C*}$ is an exact
infinitesimal Noether symmetry and consequently $\widetilde{X}$ is a natural
infinitesimal symmetry. Moreover, the function $\theta (\widetilde{X})$ is a
conservation law for the hamiltonian $H$.
\end{theorem}

\textbf{Proof}. We have that
$\mathcal{L}_{\widetilde{X}^{C*}}H=\widetilde{X} ^{C*}(H)=0$.
Next, we prove that $\mathcal{L}_{\widetilde{X}^{C*}}\theta =0$
using the computation in local coordinates.
\begin{eqnarray*}
(\mathcal{L}_{\widetilde{X}^{C*}}\theta )\left( \frac \partial {\partial
x^i}\right) &=&\widetilde{X}^{C*}\left( \theta \left( \frac \partial
{\partial x^i}\right) \right) -\theta \left[ \widetilde{X}^{C*},\frac
\partial {\partial x^i}\right] = \\
&=&\widetilde{X}^{C*}\left( p_i\right) -\theta \left(
-\frac{\partial \widetilde{X}^j}{\partial x^i}\frac \partial
{\partial x^j}+\frac{\partial ^2 \widetilde{X}^j}{\partial
x^k\partial x^i}p_j\frac \partial {\partial
p_k}\right) \\
&=&-\frac{\partial \widetilde{X}^j}{\partial x^i}p_j+\frac{\partial
\widetilde{X}^j}{\partial x^i}p_j=0,
\end{eqnarray*}
\begin{equation*}
(\mathcal{L}_{\widetilde{X}^{C*}}\theta )\left( \frac \partial
{\partial p_i}\right) =\widetilde{X}^{C*}\left( \theta \left(
\frac \partial {\partial p_i}\right) \right) -\theta \left[
\widetilde{X}^{C*},\frac \partial {\partial p_i}\right] =-\theta
\left( \frac{\partial \widetilde{X}^i}{
\partial x^j}\frac \partial {\partial p_j}\right) =0.
\end{equation*}
It results that $0=d\mathcal{L}_{\widetilde{X}^{C*}}\theta
=\mathcal{L}_{ \widetilde{X}^{C*}}d\theta
=\mathcal{L}_{\widetilde{X}^{C*}}\omega $ and $
\widetilde{X}^{C*}$ is an exact infinitesimal Noether symmetry.
Using Proposition 3.5 we have that $\widetilde{X}^{C*}$ is an
infinitesimal symmetry and consequently, $\widetilde{X}$ is a
natural infinitesimal symmetry. Moreover, according to Theorem 3.8
for $f=0$ it results that
\begin{equation*}
\theta \left( \widetilde{X}^{C*}\right) =p_i\widetilde{X}^i,
\end{equation*}
is a conservation law for the Hamiltonian $H$.

\subsection{\textbf{Example}}

Let us consider the following distributional system in $\Bbb{R}^2$
(driftless control affine system):
\begin{equation*}
\left\{
\begin{array}{l}
\dot x^1=u^1+u^2x^1 \\
\dot x^2=u^2
\end{array}
\right.
\end{equation*}
Let $x_0$ and $x_1$ be two points in $\Bbb{R}^2$. An optimal
control problem consists of finding the trajectories of our
control system which connect $ x_0 $ and $x_1$ and minimizing the
Lagrangian
\begin{equation*}
{\min }\int_0^TL(u(t))dt,\ L(u)=\frac 12\left( (u^1)^2+(u^2)^2\right) ,\quad
x(0)=x_0,\ x(T)=x_1,
\end{equation*}
where $\dot x^i=\frac{dx^i}{dt}$ and $u^1,u^2$ are control variables. Using
the Pontryagin Maximum Principle, we find the Hamiltonian function on the
cotangent bundle $T^{*}\Bbb{R}^2$ in the form
\begin{equation*}
H(x,p,u)=p_i\dot x^i-L=p_1\left( u^1+u^2x^1\right) +p_2u^2-\frac 12\left(
(u^1)^2+(u^2)^2\right) ,
\end{equation*}
with the condition $\frac{\partial H}{\partial u}=0,$ which leads
to $ u^1=p_1 $, $u^2=p_1x^1+p_2$. We obtain
\begin{equation*}
H(x,p)=\frac 12\left( p_1^2(1+(x^1)^2)+2p_1p_2x^1+p_2^2\right) =\frac
12\left( p_1^2+(p_1x^1+p_2)^2\right) ,
\end{equation*}
and it result
\begin{eqnarray*}
\frac{\partial H}{\partial p_1} &=&p_1(1+(x^1)^2)+p_2x^1,\quad
\frac{
\partial H}{\partial p_2}=p_1x^1+p_2, \\
\frac{\partial H}{\partial x^1} &=&p_1^2x^1+p_1p_2,\quad
\frac{\partial H}{
\partial x^2}=0.
\end{eqnarray*}
The Hamilton's equations lead to the following system of differential
equations
\begin{equation*}
\left\{
\begin{array}{l}
\dot x^1=p_1(1+(x^1)^2)+p_2x^1 \\
\dot x^2=p_1x^1+p_2, \\
\dot p_1=-p_1(p_1x^1+p_2) \\
\dot p_2=0\Rightarrow p_2=ct.
\end{array}
\right.
\end{equation*}
The Hessian matrix of $H$ with respect to $p$ is
\begin{equation*}
g^{ij}(x)=\frac{\partial ^2H}{\partial p_i\partial p_j}=\left(
\begin{array}{cc}
1+(x^1)^2 & x^1 \\
x^1 & 1
\end{array}
\right) ,\ i,j=1,2
\end{equation*}
and it results that $H$ is regular ($rank\left\| g^{ij}(x,p)\right\| =2$)
and its inverse matrix has the form
\begin{equation*}
g_{ij}(x)=\left(
\begin{array}{cc}
1 & -x^1 \\
-x^1 & 1+(x^1)^2
\end{array}
\right) .
\end{equation*}
The adapted tangent structure is given by
\begin{equation*}
\mathcal{J}_H=dx^1\otimes \frac \partial {\partial p_1}-x^1dx^1\otimes \frac
\partial {\partial p_2}-x^1dx^2\otimes \frac \partial {\partial
p_1}+(1+(x^1)^2)dx^2\otimes \frac \partial {\partial p_2}.
\end{equation*}
The $\mathcal{J}_H$-regular vector field is the Hamiltonian vector field
(16)
\begin{equation*}
\rho _H=\left( p_1(1+(x^1)^2)+p_2x^1\right) \frac \partial {\partial
x^1}+\left( p_1x^1+p_2\right) \frac \partial {\partial x^2}-\left(
p_1^2x^1+p_1p_2\right) \frac \partial {\partial p_1},
\end{equation*}
and from Proposition 3.7 is a infinitesimal Noether symmetry for the
dynamics induced by the regular Hamiltonian $H$. Moreover, if $X=p_2\frac
\partial {\partial x^2}$ then $[\rho _H,X]=0$ and it results that $X$ is an
infinitesimal symmetry for the Hamiltonian vector field $\rho _H$.

The local coefficients of the canonical nonlinear connection (17) have the
following form
\begin{eqnarray*}
\mathcal{N}_{11} &=&-(p_1x^1+p_2),\quad \mathcal{N}_{22}=-x^1\left(
p_1(1+(x^1)^2)+p_2x^1\right) , \\
\mathcal{N}_{12} &=&\mathcal{N}_{21}=x^1\left( p_1x^1+p_2\right) .
\end{eqnarray*}
By straightforward computation we obtain that
\begin{equation*}
\rho _H=\frac{\partial H}{\partial p_i}\frac \partial {\partial
x^i}+\frac{
\partial H}{\partial p_i}\mathcal{N}_{ij}\frac \partial {\partial p_j}=\frac{
\partial H}{\partial p_i}\frac \delta {\delta x^i},
\end{equation*}
and it results that the Hamiltonian vector field is a horizontal
$\mathcal{J} _H$-regular vector field. In this case we obtain that
the Jacobi endomorphism (18) is given by
\begin{equation*}
\mathcal{R}_{ij}=R_{kij}\frac{\partial H}{\partial p_k},
\end{equation*}
where $R_{kij}$ are the local coefficients from (2) of the curvature of the
canonical nonlinear connection with nonzero components
\begin{eqnarray*}
R_{121} &=&2p_1x^1+p_2=-R_{211}, \\
R_{212} &=&p_1+2p_1(x^1)^2+p_2x^1=-R_{122}.
\end{eqnarray*}
Also, $\nabla \rho _H=\mathcal{D}_{\rho _H}\rho _H=0$ and it results that
the integral curves of horizontal Hamiltonian vector field $\rho _H$ are
geodesics of the Berwald linear connection.

\textbf{Conclussions and further developments}.

The main purpose of this work is to study the symmetries of
Hamiltonian systems on the cotangent bundle using the same methods
as in the study of the symmetries for second order differential
equations on the tangent bundle. The role of the canonical tangent
structure and the semispray on $TM$ is taken by the adapted
tangent structure and the regular vector field on $ T^{*}M$, which
can be defined in the presence of a regular Hamiltonian. However,
the cotangent bundle has a canonical symplectic structure, which
can be found on the tangent bundle only in the presence of a
Lagrangian function. Also, we find the invariant equations of some
type of symmetries on $T^{*}M$ using the Jacobi endomorphism and
dynamical covarint derivative. Moreover, in the case of the
horizontal regular vector field, in particular the Hamiltonian
vector field, we prove that the dynamical covariant derivative
coincides with Berwald linear connection. It results that the
integral curves of the horizontal Hamiltonian vector field are the
geodesics of the Berwald linear connection. We find the relations
between infinitesimal symmetries, natural infinitesimal
symmetries, Newtonoid vector field, infinitesimal Noether
symmetries and conservation laws on $T^{*}M$ and show when one of
them will imply the others. In the last part of the paper an
examples from optimal control theory is given. As further
developments, we can use the dynamical covariant derivative and
Jacobi endomorphism in the study of symmetries for $k$-symplectic
Hamiltonian systems.

Author's address:

University of Craiova,

Dept. of Statistics and Economic Informatics

13, Al. I. Cuza, st. Craiova 200585, Romania

e-mail: liviupopescu@central.ucv.ro; liviunew@yahoo.com

\

\end{document}